\documentclass[12pt]{article}
\usepackage{amsmath,amscd,amsthm,amssymb}

\def\rn{\R^n}

\def\R{{\mathbb{R}}}
\def\Rn{{\mathbb{R}^n}}
\def\A{\mathbb A}

\def\a {\alpha}
\def\i{\infty}

\def\L1loc{L_1^{\rm loc}(\Rn)}

\def\dual{\,^{^{\complement}}\!}

\newcommand{\ess}{\mathop{\rm ess \; sup}\limits}
\newcommand{\es}{\mathop{\rm ess \; inf}\limits}

\def\Rn{\mathbb R^n}
\def\R{\mathbb R}

\def\rn{\R^n}
\def\a{\alpha}

\def\i{\infty}

\def\R{\mathbb R}

\def\dual{\,^{^{\complement}}\!}

\def\R{{\mathbb{R}}}
\def\Rn{{\mathbb{R}^n}}
\def\A{\mathbb A}

\def\a {\alpha}
\def\i{\infty}

\def\L1loc{L_1^{\rm loc}(\Rn)}

\def\dual{\,^{^{\complement}}\!}

\newtheorem{thm}{Theorem}[section]
 \newtheorem{cor}[thm]{Corollary}
 \newtheorem{lem}[thm]{Lemma}
 
 \theoremstyle{definition}
 \newtheorem{defn}[thm]{Definition}
 \theoremstyle{remark}
 \newtheorem{rem}[thm]{Remark}
 
 \numberwithin{equation}{section}

\usepackage{color}

\begin{document}

\begin{center}
\Large \bf Two-type Estimates for the Boundedness of Generalized Riesz Potential Operator in the Generalized Weighted Local Morrey Spaces
\end{center}

\centerline{Abdulhamit Kucukaslan$^{a,b,}$  \footnote{
{Corresponding author.
\\
The research of Abdulhamit Kucukaslan was supported by the grant of The Scientific and Technological Research Council of Turkey, Grant TUBITAK-1059B191600675.}
\\
E-mail address: kucukaslan@pau.edu.tr (A. Kucukaslan).}}

\centerline{$^{b}$\it Institute of Mathematics of Czech Academy of Sciences,} 

\centerline{ \it 115 67, Prague, Czech Republic}

\centerline{$^{a}$\it School of Applied Sciences, Pamukkale University, 20680, Denizli, Turkey}

\

\begin{abstract}
In this paper, we prove the Spanne-type boundedness of the generalized Riesz potential operator $I_{\rho}$ from the one generalized  weighted local Morrey spaces $M^{\{x_0\}}_{p,\varphi_{1}}(w^{p}, \Rn)$ to the another one $M^{\{x_0\}}_{q,\varphi_{2}}(w^{q}, \Rn)$ with $w^{q} \in A_{1+\frac{q}{p'}}$ for $1<p< q<\infty$ and from the generalized  weighted local Morrey spaces $M^{\{x_0\}}_{1,\varphi_1}(w, \Rn)$ to the weak generalized  weighted local Morrey spaces $WM^{\{x_0\}}_{q,\varphi_2}(w^{q}, \Rn)$ with $w \in A_{1,q}$ for $1< q<\infty$.
We also prove the Adams-type boundedness of the operator $I_{\rho}$ from the weighted spaces $M_{p,\varphi^{\frac{1}{p}}}(w, \Rn)$ to the another one $M_{q,\varphi^{\frac{1}{q}}}(w, \Rn)$  with $ w \in A_{p,q}$ for $1<p<q<\infty$  and from the weighted spaces $M_{1,\varphi}(w, \Rn)$ to the weak weighted spaces $WM_{q,\varphi^{\frac{1}{q}}}(w, \Rn)$ with $w \in A_{1,q}$ for $1<q<\infty$.
\end{abstract}

\



\

\section{Introduction} 

Morrey spaces ${M}_{p,\lambda}(\Rn)$ were introduced by Morrey in \cite{Morrey} and defined as follows: For $0 \leq \lambda <n, 1 \leq p \leq \i, f \in {M}_{p,\lambda}(\Rn)$ if $f \in L^{loc}_{p}(\rn)$ and 
$$
\|f\|_{M_{p,\lambda}(\Rn)}  = \sup\limits_{x\in\Rn, r>0} r^{-\frac{\lambda}{p}}  \|f\|_{L_p(B(x,r))}<\i
$$
holds. Morrey spaces found important applications to potential theory \cite{Adams1},  elliptic equations with discountinuous coefficients \cite{Caff} and Shr\"{o}dinger equations \cite{Ruiz}.

On the other hand, on the weighted Lebesgue spaces $L_{p}(w, \Rn)$, the boundedness of some classical operators were obtained by Muckenhoupt \cite{MW}, Mukenhoupt and Wheeden \cite{Muc}, and Coifman and Fefferman \cite{CoF}. 

Recently, weighted Morrey spaces $M_{p,\kappa}(w, \Rn)$  were introduced by Komori and Shirai \cite{KoS} as follows: For $ 1 \leq p \leq \i, 0 < \kappa < 1$ and $w$ be a weight, $f \in {M}_{p,\kappa}(w, \Rn)$ if $f \in L^{loc}_{p}(w, \Rn)$ and 
$$
\|f\|_{M_{p,\kappa}(w, \Rn)}  = \sup\limits_{x\in\Rn, r>0} w(B(x,r))^{-\frac{\kappa}{p}}  \|f\|_{L_p(w, B(x,r))}<\i.
$$
They studied the boundedness of the aforementioned classical operators in these spaces. These results were extended to several other spaces in \cite{GulEMJ2012}. Weighted inequalities for fractional operators have applications to potential theory and quantum mechanics.

For a fixed $x_0 \in \Rn$ the generalized weighted local Morrey spaces \newline $M^{\{x_0\}}_{p,\varphi}(w, \Rn)$ are obtained by replacing a function $\varphi(x_0 , r)$ instead of $r^{\lambda}$ in the definition of weighted local Morrey space, which is the space of all functions $f\in L_p^{\rm loc}(w, \Rn) $ with finite norm
$$
\|f\|_{M^{\{x_0\}}_{p,\varphi}(w, \Rn)} = \sup\limits_{r>0} \varphi(x_0,r)^{-1} \, w(B(x_0,r))^{-\frac{1}{p}} \, \|f\chi_{B(x_0,r)}\|_{L_p(w, \Rn)}.
$$
During the last decades, the theory of boundedness of classical operators of the harmonic analysis in the generalized Morrey spaces $M_{p,\varphi}(\rn)$ have been  well studied by now, we refer the readers to \cite{ GulIsmKucSJFS2015,Kuc1,GKHS1,MusKuc2018} and \cite{Nakai4}.

For a measurable function $\rho : (0,\i) \to (0,\i)$ the generalized Riesz potential operator (or generalized fractional integral operator) $I_{\rho}$ is  defined by
$$
I_{\rho}f(x) = \int_{\mathbb{R}^n }\frac{\rho(|x-y|)}{|x-y|^n} \, f(y) dy
$$
for any suitable function $f$ on $\Rn$. If $\rho(t) \equiv t^{\a}$, then we get the Riesz potential operator $I_{\a}$. The generalized Riesz potential operator $I_{\rho }$ was initilally investigated in \cite{Nakai}. Nowadays many authors have been
culminating important observations about $I_{\rho }$ especially in connection with
Morrey spaces. Nakai \cite{Nakai} proved the boundedness of $I_{\rho }$ from the generalized Morrey spaces ${M}_{1,\varphi}(\Rn)$ to the spaces ${M}_{1,\psi}(\Rn)$ for suitable functions $\varphi$ and $\psi$. The boundedness of $I_{\rho }$ from the generalized Morrey spaces ${M}_{p,\varphi }(\Rn)$ to the spaces ${M}_{q,\psi}(\Rn)$ is studied by  Eridani \cite{Eridani}, Guliyev et al \cite{GulIsmKucSJFS2015}, Kucukaslan et al \cite{Kuc1, GKHS1}, Kucukaslan \cite{Kuc-2020, Kuc-2020-3}, Nakai \cite{Nakai4} and  Nakamura \cite{Nak}.

Spanne-type and Adams-type boundednesses of generalized fractional maximal operator  $M_{\rho}$ in the generalized weighted local Morrey spaces \newline 
$M^{\{x_0\}}_{p,\varphi}(w, \Rn)$  and generalized weighted Morrey spaces $M_{p,\varphi}(w, \Rn)$  were studied in \cite{Kuc-2020-3}.  But, Spanne-type and Adams-type boundedness of the generalized Riesz potential operator $I_{\rho}$ in the spaces $M^{\{x_0\}}_{p,\varphi}(w, \Rn)$ and $M_{p,\varphi}(w, \Rn)$ have not been studied, yet.

Spanne \cite{Peetre} and Adams \cite{Adams1} studied boundedness of the Riesz potential in Morrey spaces.
Their results, can be summarized as follows.

{\bf{Theorem A.}}  {\textit{(Spanne, but published by Peetre \cite{Peetre}) \label{Peetre1}
		Let $0<\alpha<n$, $1<p<\frac{n}{\alpha}$, $0<\lambda<n-\alpha p$. Moreover, let
		$\frac{1}{p}-\frac{1}{q}=\frac{\alpha}{n}$ and $\frac{\lambda}{p}=\frac{\mu}{q}$.
		Then for $p>1$, the operator $I_{\alpha}$ is bounded from $M_{p,\lambda}(\Rn)$
		to $M_{q,\mu}(\Rn)$ and for $p=1$, $I_{\alpha}$ is bounded from $M_{1,\lambda}(\Rn)$  to $WM_{q,\mu}(\Rn)$.}}

{\bf{Theorem B.}}  {\textit{ (Adams \cite{Adams1}) \label{Adams1}
		Let $0<\alpha<n$, $1<p<\frac{n}{\alpha}$, $0<\lambda<n-\alpha p$ and
		$\frac{1}{p}-\frac{1}{q}=\frac{\alpha}{n-\lambda}$.
		Then for $p>1$, the operator $I_{\alpha}$ is bounded from $M_{p,\lambda}(\Rn)$
		to $M_{q,\lambda}(\Rn)$ and for $p=1$, $I_{\alpha}$ is bounded from $M_{1,\lambda}(\Rn)$  to $WM_{q,\lambda}(\Rn)$.
}}

In the following theorems which were proved in \cite{GulIsmKucSJFS2015}, we give Spanne and Adams type results for the boundedness of operator $I_{\rho}$ on the generalized   local Morrey spaces $M^{\{x_0\}}_{p,\varphi}(\Rn)$ and on the generalized Morrey spaces $M_{p,\varphi}(\Rn)$, respectively.

{\bf{Theorem C.}}  {\textit{(Spanne type result \cite{GulIsmKucSJFS2015}) \label{3.4.Pot}
		Let $x_0 \in \Rn$, $1 \le p <q< \infty$, the function $\rho$ satisfy the conditions \eqref{badf11}, \eqref{ank11} and [\eqref{prag38} for $w=1$]. Let also $(\varphi_1,\varphi_2)$ satisfy the conditions
		\begin{equation*}\label{VZACS}
		\es_{t<s<\infty} \varphi_1(x_0,s) s^{\frac{n}{p}}  \le  C \,\varphi_2\big(x_0,\frac{t}{2}\big) \, t^{\frac{n}{q}},
		\end{equation*}
		\begin{equation*}\label{eq3.6.VZPot}
		\int_{r}^{\infty} \Big(\es_{t<s<\infty} \varphi_1(x_0,s) s^{\frac{n}{p}} \Big)  \frac{\rho(t)}{t^{\frac{n}{p}+1}} dt \le  C \,\varphi_2(x_0,r),
		\end{equation*}
		where $C$ does not depend on $x_0$ and $r$.
		Then the operator $I_{\rho}$ is bounded from one generalized local Morrey spaces $M_{p,\varphi_1}^{\{x_0\}}(\Rn)$ to another one $M_{q,\varphi_2}^{\{x_0\}}(\Rn)$ for $p>1$
		and from the spaces $M_{1,\varphi_1}^{\{x_0\}}(\Rn)$ to the weak space $WM_{q,\varphi_2}^{\{x_0\}}(\Rn)$ for $p=1$.
}}

{\bf{Theorem D.}}  {\textit{(Adams type result \cite{GulIsmKucSJFS2015})  \label{Ris1}
		Let $1 \le p < \infty$, $q>p$, $\rho(t)$ satisfy the conditions \eqref{ank11} and [\eqref{prag38} for $w=1$]. Let also $\varphi(x,t)$ satisfy the conditions 
		\begin{equation*}\label{eq3.6.VZMaxZ}
		\sup_{r<t<\infty} \varphi(x,t) \le C \, \varphi(x,r),
		\end{equation*}
		\begin{equation*}\label{eq3.6.VXX}
		\int_{r}^{\infty} \varphi(x,t)^{\frac{1}{p}} \, \frac{\rho(t)}{t} dt \le C \rho(r)^{-\frac{p}{q-p}},
		\end{equation*}
		where $C$ does not depend on $x \in \Rn$ and $r>0$. Then the operator $I_{\rho}$ is bounded from the one generalized  Morrey space  $M_{p,\varphi^{\frac{1}{p}}}(\Rn)$ to another one
		$M_{q,\varphi^{\frac{1}{q}}}(\Rn)$ for $p>1$ and from the space $M_{1,\varphi}(\Rn)$ to the weak space $WM_{q,\varphi^{\frac{1}{q}}}(\Rn)$ for $p=1$.
}}

In this study, by using the method given in \cite{GulDoc},  we prove the Spanne and Adams type estimates for the boundedness of generalized Riesz potential operator $I_{\rho}$ on the generalized weighted local Morrey spaces $M^{\{x_0\}}_{p,\varphi}(w^{p}, \Rn)$ with $1\leq p<q<\infty,w^{q} \in A_{1+\frac{q}{p'}}$ belonging to Muckenhoupt-Weeden class $A_{p,q}$.
We find conditions on the triple $(\varphi_1,\varphi_2,\rho)$ which ensure the Spanne-type boundedness of the operator $I_{\rho}$ from one generalized  weighted local Morrey spaces $M^{\{x_0\}}_{p,\varphi_{1}}(w^{p}, \Rn)$ to another $M^{\{x_0\}}_{q,\varphi_{2}}(w^{q}, \Rn)$ with $w^{q} \in A_{1+\frac{q}{p'}}$ for $1<p< q<\infty$ and from $M^{\{x_0\}}_{1,\varphi_1}(w, \Rn)$ to the weighted weak space $WM^{\{x_0\}}_{q,\varphi_2}(w^{q}, \Rn)$ with $w \in A_{1,q}$ for $1< q<\infty$ (see Theorem \ref{kucWEG}). We also find conditions on the pair $(\varphi, \rho)$ which ensure the Adams-type boundedness of $I_{\rho}$ from $M_{p,\varphi^{\frac{1}{p}}}(w, \Rn)$ to $M_{q,\varphi^{\frac{1}{q}}}(w, \Rn)$ for $1<p<q<\infty,w \in A_{1+\frac{q}{p'}}$ and from $M_{1,\varphi}(w, \Rn)$ to $WM_{q,\varphi^{\frac{1}{q}}}(w, \Rn)$ for $1<q<\infty, w \in A_{1,q}$ (see Theorem \ref{kuc22}).

In all cases the conditions for the boundedness of $I_{\rho}$ are given in terms of Zygmund-type
integral inequalities on the all $\varphi$ functions and $r$
which do not assume any assumption on monotonicity of
$\varphi_1(x,r)$, $\varphi_2(x,r)$ and
$\varphi(x,r)$ in $r$.

By $A \lesssim B$ we mean that $A \le C B$ with some positive constant $C$ independent of appropriate quantities. If $A \lesssim B$ and $B \lesssim A$, we
write $A\approx B$ and say that $A$ and $B$ are  equivalent.
\section{Preliminaries}
Let $x \in \mathbb{R}^n $ and $r > 0$, then we denote by $B(x,r)$ the open ball centered at $x$ of radius $r$,
and by ${\,^{^{\complement}}\!} B(x,r)$ denote its complement. Let $|B(x,r)|$ be the Lebesgue measure of the ball $B(x,r)$. A weight function is a locally integrable function on $\Rn$ which takes
values in $(0,\i)$ almost everywhere. For a weight $w$ and a measurable set
$E$, we define $w(E) =\int_{E} w(x)dx,$ the Lebesgue measure of $E$ by $|E|$ and the characteristic function of $E$ by $\chi_E$.
If $w$ is a weight function, for all $f \in L_1^{loc}$ we denote by $ L_p^{loc}(w, \Rn)$ the weighted Lebesgue space defined by the norm
$$
\|f\chi_{B(x,r)}\|_{L_{p}(w, \Rn)} = \left(\int_{B(x,r)}|f(x)|^{p}w(x)dx\right)^{\frac{1}{p}}<\i,
$$
when $1 \leq p < \i$ and by
$$
\|f\chi_{B(x,r)}\|_{L_{\i}(w, \Rn)} = \ess_{x \in B(y,r)} \left|f(x)w(x)\right|<\i,
$$
when $p=\i$.

We recall that  a weight function $w$ belongs to the Muckenhoupt-Wheeden class $A_{p,q}$ (see \cite{{MW}}) for
$1 < p< q < \i$, if
\begin{align*}
\sup_{B}\left(\frac{1}{\left|B\right|}\int_{B} w(x)^q dx\right)^\frac{1}{q}\left(\frac{1}{\left|B\right|}\int_{B} w(x)^{-p'} dx\right)^\frac{1}{p'}\leq C,
\end{align*}
where the supremum is taken with respect to all balls $B$ and $C>0$. Note that, for all balls $B$ by H\"{o}lder's inequality we get
\begin{align}\label{kuc2}
|B|^{\frac{^1}{p}-\frac{1}{q}-1} \|w\|_{L_{q}(B)}\|w^{-1}\|_{L_{p'}(B)} \geq 1.
\end{align}
If $p=1$, $w$ is in the $A_{1,q}$ with $1<q<\i$ if
\begin{align*}
\sup_{B}\left(\frac{1}{\left|B\right|}\int_{B} w(x)^q dx\right)^\frac{1}{q}\left({\ess}_{x \in B}\frac{1}{w(x)}\right)\leq C,
\end{align*}
where the supremum is taken with respect to all balls $B$ and $C>0$.

The weight function $w$ satisfies the reverse doubling condition if there exist constants $\alpha_{1}>1$ and $\alpha_{2}<1$ such that 
\begin{align}\label{kuc44}
w(B(x,r))\leq \alpha_{2} w(B(x, \alpha_{1} r) )
\end{align}
for arbitrary $x \in \Rn$ and $r>0$.
\begin{lem}\cite{{GC},{Graf}}
	If $w \in A_{p,q}$ with $1 < p<q < \i$, then the following statements are true.
	\\
	$(i)$ $w^{q} \in A_{r}$ with $r=1+\frac{q}{p'}$.
	\\
	$(ii)$ $w^{-p'} \in A_{r'}$ with $r'=1+\frac{p'}{q}$.
	\\
	$(iii)$ $w^{p} \in A_{s}$ with $s=1+\frac{p}{q'}$.
	\\
	$(iv)$ $w^{-q'} \in A_{s'}$ with $s'=1+\frac{q'}{p}$.
\end{lem}

We find it convenient to define the generalized weighted local Morrey spaces in
the form as follows.
\begin{defn}
	Let $1 \le p < \infty$ and $\varphi(x,r)$ be a positive measurable function on $\Rn \times (0,\infty)$. For any fixed $x_0 \in \Rn$ we denote by
	$M_{p,\varphi}^{\{x_0\}}(w,\Rn)$ the generalized weighted local Morrey space, the space of all
	functions $f\in L_p^{\rm loc}(w, \Rn)$ with finite quasinorm
	$$
	\|f\|_{M_{p,\varphi}^{\{x_0\}}(w, \Rn)} = \|f(x_0+\cdot)\|_{M_{p,\varphi}(w, \Rn)}.
	$$
	Also by $ WM_{p,\varphi}^{\{x_0\}}(w, \Rn)$ we denote the weak
	generalized weighted local Morrey space of all functions $f\in WL_p^{\rm loc}(w, \Rn)$ for which
	$$
	\|f\|_{WM_{p,\varphi}^{\{x_0\}}(w, \Rn)}  = \|f(x_0+\cdot)\|_{WM_{p,\varphi}(w, \Rn)} < \infty.
	$$
\end{defn}

According to this definition, we recover the weighted local Morrey space $M_{p,\lambda}^{\{x_0\}}(w, \Rn)$ and weighted weak local Morrey space $WM_{p,\lambda}^{\{x_0\}}(w, \Rn)$ under the choice
$\varphi(x_0,r)=r^{\frac{\lambda-n}{p}}$:
$$
M_{p,\lambda}^{\{x_0\}}(w, \Rn)=M_{p,\varphi}^{\{x_0\}}(w, \Rn)\Big|_{\varphi(x_0,r)=r^{\frac{\lambda-n}{p}}},
$$
$$
WM_{p,\lambda}^{\{x_0\}}(w, \Rn)=WM_{p,\varphi}^{\{x_0\}}(w, \Rn)\Big|_{\varphi(x_0,r)=r^{\frac{\lambda-n}{p}}}.
$$

We denote by $L_{\infty}(w,(0,\infty))$ the space of all
functions $g(t)$, $t>0$ with finite norm
$$
\|g\|_{L_{\infty}(w,(0,\infty))} = \sup_{t>0}w(t)g(t)
$$
and $L_{\infty}(0,\infty) \equiv L_{\infty}(1,(0,\infty))$.
Let ${\mathfrak S}(0,\i)$ be the set of all Lebesgue-measurable
functions on $(0,\i)$ and ${\mathfrak S}^+(0,\i)$ its subset
consisting of all nonnegative functions on $(0,\i)$. We denote by
${\mathfrak S}^+\!(0,\i;\!\uparrow\!)\!$ the cone of all functions in
${\mathfrak S}^+(0,\i)$ which are non-decreasing on $(0,\i)$ and
$$
\A=\left\{\varphi \in {\mathfrak S}^+(0,\i;\uparrow):
\lim_{t\rightarrow 0+}\varphi(t)=0\right\}.
$$
The following theorem was proved in \cite{GulIsmKucSJFS2015} which we will use while proving our main results.
\begin{thm}\label{thm5.1XM}
	Let $w_1$, $w_2$ be non-negative measurable functions satisfying $0<\|w_1\|_{L_{\i}(t,\i)}<\i$ for any $t>0$. Then the identity operator $I$ is bounded from $L_{\i}(w_1,(0,\i))$ to $L_{\i}(w_2,(0,\i))$ on the cone $\A$ if and only if
	\begin{equation*}\label{gfa01}
	\begin{split}
	\left\|w_2 \left( \| w_1 \|^{-1}_{L_{\i}(\cdot,\i)}\right)\right\|_{L_{\i}(0,\i)}<\i.
	\end{split}
	\end{equation*}
\end{thm}

We will  use the following statement on the boundedness of the weighted Hardy operator
$$
H_{w}g(t):=\int_t^{\i} g(s) w(s) d\mu(s),~ \ \  0<t<\infty,
$$
where $w$ is weight and $d\mu(s)$ is a non-negative Borel measure on $(0,\i)$.

\
The following theorem was proved in \cite{CPSS}.
\begin{thm}\label{thm3.2.}
	Let $w_1$, $w_2$ and $w$ be weights on  $(0,\infty)$ and $w_1(t)$ be bounded outside a neighborhood of the
	origin. The inequality
	\begin{equation} \label{vav01}
	\ess_{t>0} w_2(t) H_{w} g(t) \leq C \ess_{t>0} w_1(t) g(t)
	\end{equation}
	holds for some $C>0$ for all non-negative and non-decreasing $g$ on $(0,\i)$ if and
	only if
	\begin{equation} \label{vav02}
	B:= \sup _{t>0} w_2(t)\int_t^{\infty} \frac{w(s) ds}{\ess_{s<\tau<\infty} w_1(\tau)}<\infty.
	\end{equation}
	Moreover, the value  $C=B$ is the best constant for  \eqref{vav01}.
\end{thm}

\begin{rem}\label{rem2.3.}
	In \eqref{vav01} and \eqref{vav02} it is assumed that $\frac{1}{\i}=0$ and $0 \cdot \i=0$.
\end{rem}

\section{ Spanne-type result for the operator $I_{\rho}$ in the spaces $M_{p,\varphi}^{\{x_0\}}(w^{p}, \Rn)$}

We assume that
\begin{equation} \label{badf11}
\int_{1}^{\i} \frac{\rho(t)}{t^n} \frac{dt}{t} < \infty,
\end{equation}
so that the generalized Riesz potential $I_{\rho}f$ is well defined,
at least for characteristic functions $1/|x|^{2n}$ of complementary balls:
$$
f(x)=\frac{\chi_{\Rn \setminus B(0,1)}(x)}{|x|^{2n}}.
$$
In addition, we shall also
assume that $\rho$ satisfies the growth condition: there exist constants $C > 0$ and $0<2k_1<k_2<\i$ such that
\begin{equation} \label{ank11}
\sup\limits_{r < s \le 2r} \frac{\rho(s)}{s^n} \le C \int_{k_1 r}^{k_2 r} \frac{\rho(t)}{t^n} \frac{dt}{t}, ~~ r>0.
\end{equation}

This condition is weaker than the usual doubling condition for the function $\frac{\rho(t)}{t^n}$ : there exists a constant $C> 0$ such that
\begin{equation} \label{ank21}
\frac{1}{C} \frac{\rho(t)}{t^n} \le \frac{\rho(r)}{r^n} \le C \frac{\rho(t)}{t^n},
\end{equation}
whenever $r$ and $t$ satisfy $r$, $t>0$ and $\frac{1}{2} \le \frac{r}{t} \le 2$.

In the sequel for the generalized  Riesz potential operator $I_{\rho}$ we always
assume that $\rho$ satisfies the conditions \eqref{ank11} and, then denote
the set of all such functions by $\widetilde{G}_{0}$. We will write, when $\rho \in \widetilde{G}_{0},$
\[
\widetilde{\rho}(r):= C r^{n} \int_{r}^{\infty} \frac{\rho(t)}{t^{n}} \, \frac{dt}{t}.
\]

The following lemma is valid for the operator  $I_{\rho}$.
\begin{lem}\cite{GGK}\label{prag37}
	Let $w^{q} \in A_{1+\frac{q}{p'}}$ satisfies \eqref{kuc44}, the function $\rho$ satisfies the conditions\eqref{badf11}-\eqref{ank21}, 
	and $f\in L_{1}^{loc} (w, \Rn)$. Then there exist $C>0$ for all $B(x,r)\subset \rn$  such that the inequality
	\begin{equation}\label{prag38}
	\sup_{x \in \Rn, r>0}\frac{\rho (r)}{r^n} \left(\int_{B(x,r)}w^{q}(x)dx\right)^{\frac{1}{q}}
	\left(\int_{B(x,r)}w(x)^{-p'}dx\right)^{\frac{1}{p'}}\leq C
	\end{equation}
	is necessary and sufficient condition for the boundedness of generalized Riesz potential operator $I_{\rho}$  from $L_{p}(w^{p}, \Rn)$ to $WL_{q}(w^{q}, \Rn)$ for  $1 \leq  p < q < \i$, and from $L_{p}(w^{p}, \Rn)$ to $L_{q}(w^{q}, \Rn)$ for  $1 <  p < q < \i, w^{q} \in A_{1+\frac{q}{p'}}$, where the constant $C$ does not depend on $f$.
\end{lem}

The following is weighted local $L_{p}(\Rn)$-estimate for the operator $I_{\rho}$.
\begin{lem}\label{lem3.3.Pot}
	Let fixed $x_{0} \in \Rn$, and $1 \leq  p < q < \i$, $w \in A_{1+\frac{q}{p'}}$ and $\rho(t)$ satisfy the conditions \eqref{badf11} and \eqref{ank11}.
	
	If the condition \eqref{prag38} is fulfill, then the inequality
	\begin{align}\label{eq3.5.}
	&\|I_{\rho} f\chi_{B(x_{0},r)}\|_{L_q(w^{q}, \Rn)} \lesssim   \|f\chi_{B(x_{0},2r)}\|_{L_p(w^{p}, \Rn)}\notag
	\\
	&+ \left(w^{q}(B(x_{0},r))\right)^{\frac{1}{q}} \int_{2r}^{\i} \|f\chi_{B(x_{0},t)}\|_{L_p(w^{p}, \Rn)} \left(w^{q}(B(x_{0},r))\right)^{-\frac{1}{q}} \frac{\rho(t)}{t^{n}} \frac{dt}{t}
	\end{align}
	holds for the ball $B(x_{0},r)$ and for all $f\in L_{p}^{loc}(\Rn,w)$.
	
	If the condition \eqref{prag38} is fulfill, then for $p=1$ the inequality
	\begin{align}\label{eq3.5.WXPot}
	&\|I_{\rho} f\chi_{B(x_{0},r)}\|_{WL_q(w^{q})} \lesssim  \|f\chi_{B(x_{0},2r)}\|_{L_1(w)} \notag
\\
&	+\left(w^{q}(B(x_{0},r))\right)^{\frac{1}{q}} \int_{2r}^{\i} \|f\chi_{B(x_{0},t)}\|_{L_1(w)} \, \left(w^{q}(B(x_{0},t))\right)^{-\frac{1}{q}} \frac{\rho(t)}{t^{n}} \frac{dt}{t}
	\end{align}
	holds for the ball $B(x_{0},r)$ and for all $f \in L_{1}^{loc}(w)$.
\end{lem}
\begin{proof}
	Let $1 \leq  p < q < \i$ and $w \in A_{1+\frac{q}{p'}}$.
	For fixed $x_{0} \in\Rn$, set $B\equiv B(x_{0},r)$ for
	the ball centered at $x_{0}$ and of radius $r$.
	Write $f=f_1+f_2$ with $f_1=f\chi_{2B}$ and $f_2=f\chi_{\dual {(2B)}}$. Hence, by the Minkowski inequality we have
	$$
	\|I_{\rho} f\chi_{B}\|_{L_{q}(w^{q}, \Rn)} \le \|I_{\rho} f_1\chi_{B}\|_{L_{q}(w^{q}, \Rn)} +\|I_{\rho} f_2\chi_{B}\|_{L_{q}(w^{q}, \Rn)}.
	$$
	
	Since $f_1\in L_p(w^{p}, \Rn)$, $I_{\rho} f_1\in L_q(w^{q}, \Rn)$  and from condition \eqref{prag38} we get the boundedness of $I_{\rho}$ from  $L_p(w^{p}, \Rn)$ to $L_q(w^{q}, \Rn)$ (see Lemma \ref{prag37}) and it follows that:
	\begin{equation*}
	\|I_{\rho} f_1\chi_{B}\|_{L_q(w^{q}, \Rn)}\leq \|I_{\rho} f_1\|_{L_q(w^{q}, \Rn)}\leq
	C\|f_1\|_{L_p( w^{p}, \Rn)}=C\|f\chi_{2B}\|_{L_p(w^{p}, \Rn)},
	\end{equation*}
	where constant $C >0$ is independent of $f$.
	
	It's clear that $z\in B$, $y\in \dual {(2B)}$ implies $\frac{1}{2}|x_{0}-y|\le|z-y|\le \frac{3}{2}|x_{0}-y|$.
	Then from conditions \eqref{badf11}, \eqref{ank11} and by Fubini's theorem we have
	\begin{equation*}
	\begin{split}
	|I_{\rho}f_{2}(z)|
	& \lesssim \int_{\dual {(2B)}} \frac{\rho(|x-y|)}{|x-y|^{n}} |f(y)| dy
		\lesssim \int_{2r}^{\i}\int_{B(x_{0},t) }|f(y)|dy \frac{\rho(t)}{t^{n}} \frac{dt}{t}.
	\end{split}
	\end{equation*}
	Applying H\"older's inequality and from \eqref{kuc2}, we get
	\begin{align*}\label{kuc21}
	\int_{\dual {(2B)}} &\frac{\rho(|x-y|)}{|x-y|^{n}} |f(y)| dy\notag
	\\
	&\lesssim
	\int_{2r}^{\i}\|f\chi_{B(x_{0},t)}\|_{L_p(w^{p}, \Rn)} \|w^{-1}\chi_{B(x_{0},t)}\|_{L_{p'}( \Rn)}\frac{\rho(t)}{t^{n}} \frac{dt}{t}\notag
	\\
	&\lesssim
	\int_{2r}^{\i}\|f\chi_{B(x_{0},t)}\|_{L_p(w^{p}, \Rn)} \left(w^{q}(B(x_{0},t))\right)^{-\frac{1}{q}}\frac{\rho(t)}{t^{n}} \frac{dt}{t}.
	\end{align*}
	Moreover, for all $p\in [1,\i)$ the inequality
	\begin{align} 
	&\|I_{\rho} f_2\chi_{B}\|_{L_q(w^{p}, \Rn)}\notag
	\\
	&\lesssim \left(w^{q}(B)\right)^{\frac{1}{q}}\int_{2r}^{\i}\|f\chi_{B(x_{0},t)}\|_{L_p(w^{p}, \Rn)} \left(w^{q}(B(x_{0},t))\right)^{-\frac{1}{q}} \frac{\rho(t)}{t^{n}} \frac{dt}{t}
	\end{align}
	is valid. Thus
	\begin{align*}
	\|I_{\rho} f\chi_{B}\|&_{L_q(w^{q}, \Rn)}
	\lesssim \|f\chi_{2B}\|_{L_p(w^{p}, \Rn)}
	\\
&	+\left(w^{q}(B)\right)^{\frac{1}{q}}\int_{2r}^{\i}\|f\chi_{B(x_{0},t)}\|_{L_p(w^{p}, \Rn)} \|w^{-1}\chi_{B(x_{0},t)}\|_{L_{p'}(\Rn)} \frac{\rho(t)}{t^{n}} \frac{dt}{t}.
	\end{align*}
	On the other hand,
\begin{align}\label{ves2Pot}
	\|f\chi_{2B}\|&_{L_p(w^{p}, \Rn)}\notag
	\approx
	\frac{r^{\frac{n}{p}}}{\rho(r)}\|f\chi_{2B}\|_{L_p(w^{p}, \Rn)}\int_{r}^{\i} \frac{\rho(t)}{t^{n}} \frac{dt}{t}\notag
	\\
	&\leq \frac{r^{\frac{n}{p}}}{\rho(r)}\int_{2r}^{\i} \|f\chi_{2B(x_{0},t)}\|_{L_p(w^{p}, \Rn)} \frac{\rho(t)}{t^{n}} \frac{dt}{t}\notag
	\\
	&\lesssim \left(w^{q}(B)\right)^{\frac{1}{q}} \|w^{-1}\chi_{B}\|_{L_{p'}( \Rn)} \int_{2r}^{\i} \|f\chi_{2B(x_{0},t)}\|_{L_p(w^{p})} \frac{\rho(t)}{t^{n}} \frac{dt}{t}\notag
	\\
	&\lesssim \left(w^{q}(B)\right)^{\frac{1}{q}} \int_{2r}^{\i} \|f\chi_{2B(x_{0},t)}\|_{L_p(w^{p}, \Rn)} \|w^{-1}\chi_{B(x_{0},t)}\|_{L_{p'}( \Rn)} \frac{\rho(t)}{t^{n}} \frac{dt}{t}\notag
	\\
	&\lesssim  \left(w^{q}(B)\right)^{\frac{1}{q}} \int_{2r}^{\i} \|f\chi_{2B(x_{0},t)}\|_{L_p(w^{p}, \Rn)} \left(w^{q}(B(x_{0},t))\right)^{-\frac{1}{q}} \frac{\rho(t)}{t^{n}} \frac{dt}{t}.
\end{align}
	Hence by
	\begin{equation*} 
	\|I_{\rho} f\chi_{B}\|_{L_q(w^{p}, \Rn)}\lesssim
	\left(w^{q}(B)\right)^{\frac{1}{q}}\int_{2r}^{\i}\|f\chi_{B(x_{0},t)}\|_{L_p(w^{p}, \Rn)} \left(w^{q}(B(x_{0},t))\right)^{-\frac{1}{q}} \frac{\rho(t)}{t^{n}} \frac{dt}{t}
	\end{equation*}
	we get the inequality \eqref{eq3.5.}. 
	
	Now let $p=1$ and $w \in A_{1,q}$. In this case by \eqref{prag38} we obtain
	\begin{align} \label{gfr9Pot}
	\|I&_{\rho} f_1\chi_{B}\|_{WL_q(w^{q},\Rn)}  \leq \|I_{\rho} f_1\|_{WL_q(w^{q},\Rn)}
	\lesssim \|f_1\|_{L_1(w,\Rn)} = \|f\chi_{B}\|_{L_1(w,\Rn)}\notag
	\\
	&\approx\frac{r^{n}}{\rho(r)}\|f\chi_{2B}\|_{L_1(w,\Rn)}\int_{r}^{\i} \frac{\rho(t)}{t^{n}} \frac{dt}{t}\notag
	\\
	&\leq \frac{r^{n}}{\rho(r)}\int_{2r}^{\i} \|f\chi_{2B(x_{0},t)}\|_{L_1(w,\Rn)} \frac{\rho(t)}{t^{n}} \frac{dt}{t}\notag
	\\
	&\lesssim \left(w^{q}(B)\right)^{\frac{1}{q}} \|w^{-1}\chi_{B}\|_{L_{\i}} \int_{2r}^{\i} \|f\chi_{2B(x_{0},t)}\|_{L_1(w,\Rn)} \frac{\rho(t)}{t^{n}} \frac{dt}{t}\notag
	\\
	&\lesssim \left(w^{q}(B)\right)^{\frac{1}{q}} \int_{2r}^{\i} \|f\chi_{2B(x_{0},t)}\|_{L_1(w,\Rn)} \|w^{-1}\chi_{B(x_{0},t)}\|_{L_{\i}} \frac{\rho(t)}{t^{n}} \frac{dt}{t}\notag
	\\
	&\lesssim \left(w^{q}(B)\right)^{\frac{1}{q}} \int_{2r}^{\i} \|f\chi_{2B(x_{0},t)}\|_{L_1(w,\Rn)} \left(w^{q}(B(x_{0},t))\right)^{-\frac{1}{q}} \frac{\rho(t)}{t^{n}} \frac{dt}{t}.
	\end{align}
	
	Then from \eqref{ves2Pot} and \eqref{gfr9Pot} we get the inequality \eqref{eq3.5.WXPot}.
\end{proof}

The following theorem one of the main result of our paper, in which we prove the Spanne-type estimate for the boundedness of generalized fractional integral operator $I_{\rho}$  from generalized weighted local Morrey spaces $M_{p,\varphi_1}^{\{x_0\}}(w^{p},\Rn)$ to $M_{q,\varphi_2}^{\{x_0\}}(w^{q},\Rn)$.

\begin{thm}\label{kucWEG}
	Let fixed $x_0 \in \Rn$, $1 \le p <q< \infty$, $w \in  A_{1+\frac{q}{p'}}$, the function $\rho$ satisfy the conditions \eqref{badf11}, \eqref{ank11} and \eqref{prag38}. Let also $(\varphi_1,\varphi_2)$ satisfy the conditions
	\begin{equation}\label{kuc13}
	\es_{t<s<\infty} \varphi_1(x_0,s) s^{\frac{n}{p}}  \le  C \,\varphi_2\big(x_0,\frac{t}{2}\big) \, t^{\frac{n}{q}},
	\end{equation}
	\begin{equation}\label{kuc14}
	\int_{r}^{\infty} \frac{\es_{t<s<\infty} \varphi_1(x_0,s)(w^{p}(B(x_0,s)))^{\frac{1}{p}}\rho(t)}{(w^{q}(B(x_0,t)))^{\frac{1}{q}}t^{\frac{n}{p}}}\frac{dt}{t} \le  C \,\varphi_2(x_0,r),
	\end{equation}
	where $C$ does not depend on $x$ and $r$.
	Then the operator $I_{\rho}$ is bounded from $M_{p,\varphi_1}^{\{x_0\}}(w^{p},\Rn)$ to $M_{q,\varphi_2}^{\{x_0\}}(w^{q},\Rn)$ for $p>1$
	and from $M_{1,\varphi_1}^{\{x_0\}}(w,\Rn)$ to $WM_{q,\varphi_2}^{\{x_0\}}(w^{q},\Rn)$ for $p=1$. Moreover, for $p>1$
	\begin{equation*}
	\|I_{\rho} f\|_{M_{q,\varphi_2}^{\{x_0\}}(w^{q},\Rn)} \lesssim \|f\|_{M_{p,\varphi_1}^{\{x_0\}}(w^{p},\Rn)},
	\end{equation*}
	and for $p=1$
	\begin{equation*}
	\|I_{\rho} f\|_{WM_{q,\varphi_2}^{\{x_0\}}(w^{q},\Rn)} \lesssim \|f\|_{M_{1,\varphi_1}^{\{x_0\}}(w,\Rn)}.
	\end{equation*}
\end{thm}
\begin{proof}
	Let $p>1$. Then by Theorems \ref{thm5.1XM}, \ref{thm3.2.}, Lemma \ref{lem3.3.Pot} and conditions \eqref{kuc13}-\eqref{kuc14}  we have that
	\begin{equation*}
	\begin{split}
	\|I_{\rho} f\|_{M_{q,\varphi_2}^{\{x_0\}}(w^{q},\Rn)}& =\sup_{r>0}\varphi_2(x_0,r)^{-1}(w^{q}(B(x_0,r)))^{-\frac{1}{q}} \|I_{\rho}f\|_{L_q(w^{q}, B(x_0,2r))}
	\\
	&= \sup_{r>0}\varphi_2(x_0,r)^{-1}(w^{q}(B(x_0,r)))^{-\frac{1}{q}} \|I_{\rho}f{\chi_{B(x_{0},2r)}}\|_{L_q(w^{q}, \Rn)}
	\\
	& \lesssim  \sup_{r>0}\varphi_2(x_0,r)^{-1}(w^{q}(B(x_0,r)))^{-\frac{1}{q}} \|f\|_{L_p(w^{p}, B(x_0,2r))} 
	\\
	&+
	\sup_{r>0} \varphi_2(x_0,r)^{-1} \int_r^{\i}\|f\|_{L_p(w^{p}, B(x_0,2t))} \frac{\rho(t)}{t^{\frac{n}{p}+1}} dt
	\\
	& \thickapprox  \sup_{r>0} \varphi_1(x_0,r)^{-1} (w^{p}(B(x_0,r)))^{-\frac{1}{p}} \|f\|_{L_p(w^{p}, B(x_0,r))}
	\\
	& = \|f\|_{M_{p,\varphi_1}^{\{x_0\}}(w^{p}, \Rn)}.
	\end{split}
	\end{equation*}
	Now let $p=1$ then
	\begin{equation*}
	\begin{split}
	\|I_{\rho} f\|_{WM_{q,\varphi_2}^{\{x_0\}}(w^{q}, \Rn)}& =\sup_{r>0}\varphi_2(x_0,r)^{-1}(w^{q}(B(x_0,r)))^{-\frac{1}{q}} \|I_{\rho}f\|_{WL_q(w^{q}, B(x_0,2r))}
	\\
	&= \sup_{r>0}\varphi_2(x_0,r)^{-1}(w^{q}(B(x_0,r)))^{-\frac{1}{q}} \|I_{\rho}f{\chi_{B(x_{0},2r)}}\|_{WL_q(w^{q}, \Rn)}
	\\
	& \lesssim 
	\sup_{r>0}\varphi_2(x_0,r)^{-1}(w^{q}(B(x_0,r)))^{-\frac{1}{q}} \|f\|_{L_1(w,B(x_0,2r))} 
	\\
	&+
	\sup_{r>0} \varphi_2(x_0,r)^{-1} \int_r^{\i}\|f\|_{L_1(w, B(x_0,2t))} \frac{\rho(t)}{t^{n+1}} dt
	\\
	& \thickapprox  \sup_{r>0}\varphi_1(x_0,r)^{-1}(w(B(x_0,r)))^{-1} \|f\|_{L_1(w, B(x_0,r))} 
	\\
	&= \|f\|_{M_{1,\varphi_1}^{\{x_0\}}(w, \Rn)}.
	\end{split}
	\end{equation*}
	Hence the proof is completed.
\end{proof}
\begin{cor}
	In the case $ w \equiv 1$ from Theorem \ref{kucWEG} we get Theorem C, in which we give Spanne-type result for  generalized Riesz potential operator $I_{\rho}$ on generalized local Morrey spaces $M_{p,\varphi}^{\{x_0\}}(\Rn)$ which was proved in \cite{GulIsmKucSJFS2015} (Theorem 16, p.6).
\end{cor}
\begin{cor}
	In the case $\rho(t)=t^{\a}, w \equiv 1, x \equiv x_0 $ from Theorem \ref{kucWEG} we get Spanne-type result 	for  Riesz potential operator $I_{\alpha}$ on generalized Morrey spaces $M_{p,\varphi}(\Rn)$
	which was proved in \cite{GULAKShIEOT2012} (Theorem 5.4, p.338).
\end{cor}
\begin{cor}
	In the case $\rho(t)=t^{\a}, w \equiv 1$ and $\varphi_(x_0,t)=t^\frac{{\lambda-n}}{p}$, $0<\lambda<n $ from Theorem \ref{kucWEG} we get Spanne result for  Riesz potential operator $I_{\alpha}$ on local Morrey spaces $M_{p,\lambda}^{\{x_0\}}(\Rn)$ which is variant of Theorem A proved in \cite{Peetre}.
\end{cor}
\section{ Adams-type result for the operator $I_{\rho}$ in in the spaces $M_{p,\varphi}(w)$}
The following theorem is  Adams-type estimate for generalized Riesz potential operator $I_{\rho}$ on generalized weighted Morrey spaces $M_{p,\varphi}(w, \Rn)$.

\begin{thm}\label{kuc22}
	Let $1 \le p < \infty$, $q>p$, $w \in  A_{1+\frac{q}{p'}}$, $\rho(t)$ satisfy the conditions \eqref{badf11}, \eqref{ank11} and \eqref{prag38}. Let also $\varphi(x,t)$ satisfy the conditions
	\begin{equation} \label{kuc23}
	c^{-1} \varphi(x,r)\le \varphi(x,t)\le c \, \varphi(x,r),
	\end{equation}
	\begin{equation}\label{kuc24}
	\int_{r}^{\infty} \frac{\es_{t<s<\infty} (\varphi(x,s)w(B(x,s)))^{\frac{1}{p}}}{w(B(x,s))^{\frac{1}{q}}}\frac{\rho(t)}{t^{n}}\frac{dt}{t} \leq  C \, \left( \rho(r)\right)^{-\frac{p}{q-p}},
	\end{equation}
	where $C$ does not depend on $x \in \Rn$ and $r>0$.	Then the operator $I_{\rho}$ is bounded from $M_{p,\varphi^{\frac{1}{p}}}(w, \Rn)$ to
	$M_{q,\varphi^{\frac{1}{q}}}(w, \Rn)$ for $p>1$ and from $M_{1,\varphi}(w, \Rn)$ to $WM_{q,\varphi^{\frac{1}{q}}(w, \Rn)}$ for $p=1$. Moreover, for $p>1$
	\begin{equation*}
	\|I_{\rho} f\|_{M_{q,\varphi^{\frac{1}{q}}(w, \Rn)}} \lesssim \|f\|_{M_{p,\varphi^{\frac{1}{p}}}(w, \Rn)},
	\end{equation*}
	and for $p=1$
	\begin{equation*}
	\|I_{\rho} f\|_{WM_{q,\varphi^{\frac{1}{q}}(w, \Rn)}} \lesssim \|f\|_{M_{1,\varphi}(w, \Rn)}.
	\end{equation*}
\end{thm}

\begin{proof}
	Let $1 < p < \infty, q>p, w \in  A_{1+\frac{q}{p'}}$ and $f\in M_{p,\varphi^{\frac{1}{p}}}(w, \Rn)$.
	Write $f=f_1+f_2$, where $B=B(x,r)$, $f_1=f\chi_{2B}$ and $f_2=f\chi_{{\,^{^{\complement}}\!}{(2B)}}$.
	Then we have
	$$
	I_{\rho} f(x)=I_{\rho} f_1(x)+I_{\rho} f_2(x).
	$$
	For $I_{\rho} f_1(y)$, $y \in B(x,r)$, following  Hedberg's trick (see for instance \cite{Stein93}, p. 354), we obtain
	\begin{align*}
	\left|I_{\rho} f_1(y)\right| &  \lesssim Mf(x) \rho(r),
	\end{align*}
	(see \cite{GulIsmKucSJFS2015}, for more detail).
	Thus by taking $ L_{q}(w, \Rn)-$norm we get
	\begin{align*}
	\|I_{\rho} f_1\|_{L_q(B(x,r))}(w, \Rn) & \le w(B(x,r))^{\frac{1}{q}} \left( \int_{B(y,2r)} \left(\frac{\rho(|y-z|)}{|y-z|^{n}} |f(z)| \right)^q dz \right)^{1/q}.
	\end{align*}
	
	For $I_{\rho}f_2(y)$, $y \in B(x,r)$ from \eqref{kuc2} we have
	\begin{align} \label{kirV1Zs}
	\left|I_{\rho} f_2(y)\right|  &\lesssim  \int_{{\,^{^{\complement}}\!}B(x,2r)} \frac{\rho(|y-z|)}{|y-z|^{n}} |f(z)| dz \notag
	\\
&	\lesssim 	\int_{2r}^{\i}\|f\chi_{B(x,t)}\|_{L_p(w, \Rn)} w(B(x,r))^{-\frac{1}{q}}\frac{\rho(t)}{t^{n}} \frac{dt}{t}.
	\end{align}
	Then from condition \eqref{kuc24} and inequality \eqref{kirV1Zs} for all $y \in B(x,r)$ we get
	\begin{align*} \label{kirBGZs}
	|I_{\rho} f(y)| & \lesssim \rho(r) \, Mf(x) +   	\int_{2r}^{\i}\|f\chi_{B(x,t)}\|_{L_p(w, \Rn)} w(B(x,r))^{-\frac{1}{q}} \frac{\rho(t)}{t^{n}} \frac{dt}{t}   \notag
	\\
	& \le  \rho(r) \, Mf(x) 
	\\
	& +  \|f\|_{M_{p,\varphi^{\frac{1}{p}}}(w, \Rn)} \; \int_{r}^{\infty} \varphi(x,t)^{\frac{1}{p}}w(B(x,t))^{\frac{1}{p}}{w(B(x,t))^{\frac{1}{q}}}\frac{\rho(t)}{t^{n}}\frac{dt}{t}  \notag
	\\
	& \lesssim \rho(r) \, Mf(x) + \rho(r)^{-\frac{p}{q-p}} \; \|f\|_{M_{p,\varphi^{\frac{1}{p}}}(w, \Rn)}.
	\end{align*}
	
	Hence choosing $\rho(r)=\Big(\frac{\|f\|_{M_{p,\varphi^{1/p}}(w, \Rn)}}{Mf(x)}\Big)^{\frac{q-p}{q}}$ for all $y \in B(x,r)$, we have
	$$
	|I_{\rho} f(y)| \lesssim (Mf(x))^{\frac{p}{q}} \, \|f\|_{M_{p,\varphi^{\frac{1}{p}}}(w, \Rn)}^{1-\frac{p}{q}}.
	$$
	Consequently the statement of the theorem follows in view of the  boundedness of the
	maximal operator $M$ in $M_{p,\varphi^{\frac{1}{p}}}(\Rn)$ provided in \cite{GulEMJ2012} in virtue of condition \eqref{kuc23}, hence, for $1<p<q<\infty$ we get
	\begin{align*}
&	\|I_{\rho} f\|_{M_{q,\varphi^{\frac{1}{q}}}(w, \Rn)}  = \sup_{x \in \Rn, t>0} \varphi(x,t)^{-\frac{1}{q}} w(B(x,t))^{-\frac{1}{q}} \|I_{\rho} f\|_{L_q(w, B(x,t))}
	\\
	& \lesssim \, \|f\|_{M_{p,\varphi^{\frac{1}{p}}}(w, \Rn)}^{1-\frac{p}{q}} \, \sup_{x \in \Rn, t>0} \varphi(x,t)^{-\frac{1}{q}}w(B(x,t))^{-\frac{1}{q}} \|M f\|_{L_p(w, B(x,t))}^{\frac{p}{q}}
	\\
	& = \, \|f\|_{M_{p,\varphi^{\frac{1}{p}}}(w, \Rn)}^{1-\frac{p}{q}} \,
	\left(\sup_{x \in \Rn, t>0} \varphi(x,t)^{-\frac{1}{p}} (w(B(x,t)))^{-\frac{1}{p}} \|M f\|_{L_p(w, B(x,t))}\right)^{\frac{p}{q}}
	\\
	& = \, \|f\|_{M_{p,\varphi^{\frac{1}{p}}}(w, \Rn)}^{1-\frac{p}{q}} \, \|M f\|_{M_{p,\varphi^{\frac{1}{p}}}(w, \Rn)}^{\frac{p}{q}}
	\\
	& \lesssim \, \|f\|_{M_{p,\varphi^{\frac{1}{p}}}(w, \Rn)},
	\end{align*}
	and for $1<q<\infty$ 
	\begin{align*}
	&\|I_{\rho} f\|_{WM_{q,\varphi^{\frac{1}{q}}}(w, \Rn)}  = \sup_{x \in \Rn, t>0}
	\varphi(x,t)^{-\frac{1}{q}} w(B(x,t))^{-\frac{1}{q}} \|I_{\rho} f\|_{WL_q(w, B(x,t))}
	\\
	& \lesssim \, \|f\|_{M_{1,\varphi}(w, \Rn)}^{1-\frac{1}{q}} \, \sup_{x \in \Rn, t>0}
	\varphi(x,t)^{-\frac{1}{q}} w(B(x,t))^{-\frac{1}{q}} \|M f\|_{WL_1(w, B(x,t))}^{\frac{1}{q}}
	\\
	& = \, \|f\|_{M_{1,\varphi}(w, \Rn)}^{1-\frac{1}{q}} \, \left(\sup_{x \in \Rn, t>0} \, \varphi(x,t)^{-1} w(B(x,t))^{-1} \|M f\|_{WL_1(w, B(x,t))}\right)^{\frac{1}{q}}
	\\
	& = \, \|f\|_{M_{1,\varphi}(w, \Rn)}^{1-\frac{1}{q}} \, \|M f\|_{M_{1,\varphi}(w, \Rn)}^{\frac{1}{q}}
	\\
	& \lesssim \, \|f\|_{M_{1,\varphi}(w, \Rn)}.
	\end{align*}
	Hence the proof is completed.
\end{proof}

\begin{cor}
	In the case $ w \equiv 1$ from Theorem \ref{kuc22} we get Theorem D, in which we give Adams type result	for  generalized Riesz potential operator $I_{\rho}$ on generalized  Morrey spaces $M_{p,\varphi}(\Rn)$ which was proved in \cite{GulIsmKucSJFS2015} (Theorem 22, p.7).
\end{cor}

\begin{cor}
	In the case $\rho(t)=t^{\a}, w \equiv 1, x \equiv x_0 $ from Theorem \ref{kuc22} we get Adams type result for  Riesz potential operator $I_{\alpha}$ on generalized Morrey spaces $M_{p,\varphi}(\Rn)$
	which was proved in \cite{GULAKShIEOT2012} (Theorem 5.7, p.182).
\end{cor}

\begin{cor}
	In the case $\rho(t)=t^{\a}, w \equiv 1$ and $\varphi_(x_0,t)=t^\frac{{\lambda-n}}{p}$, $0<\lambda<n $ from Theorem \ref{kuc22} we get Adams result for  Riesz potential operator $I_{\alpha}$ on local Morrey spaces $M_{p,\lambda}^{\{x_0\}}(\Rn)$ which is variant of Theorem B proved in \cite{Peetre}.
\end{cor}

\section*{Acknowledgement} 
The author would like to express his gratitude to the referees for their (his/her) very valuable comments.
\section*{Conflicts of Interest}
The author declares that there are no conflicts of interest regarding the publication of this article.

\end{document}